\documentclass[12pt,reqno]{amsart}
\usepackage{a4}
\usepackage{amssymb}
\usepackage{graphicx}
\usepackage{array}
\def\textik#1{\hbox{\footnotesize#1}}
\theoremstyle{plain}
\newtheorem{theorem}{Theorem}
\newtheorem{lemma}{Lemma}
\theoremstyle{remark}
\newtheorem{remark}{Remark}
\theoremstyle{definition}
\newtheorem{dfn}{Definition}

\def\pa#1#2{\frac{\partial#1}{\partial#2}}

\DeclareMathOperator{\rank}{rank}

\DeclareMathOperator{\Ima}{Im}
\DeclareMathOperator{\Ker}{Ker}
\DeclareMathOperator{\const}{const}

\begin{document}
\author{Rinat Kashaev, Igor Korepanov and Evgeniy Martyushev}
\title[A finite-dimensional TQFT based on group $\mathrm{PSL}(2,\mathbb C)$]{A finite-dimensional TQFT for three-manifolds based on group {\boldmath$\mathrm{PSL}(2,\mathbb C)$} and cross-ratios}
\begin{abstract}
In this paper, we begin constructing a new finite-dimensional topological quantum field theory (TQFT) for three-manifolds, based on group
$\mathrm{PSL}(2,\mathbb C)$ and its action on a complex variable by fractional-linear transformations, by providing its key ingredient --- a new type
of chain complexes. As these complexes happen to be acyclic often enough, we make use of their torsion to construct different versions of manifold
invariants. In particular, we show how to construct a large set of invariants for a manifold with boundary, analogous to the set of invariants based
on Euclidean geometric values and used in a paper by one of the authors for constructing a ``Euclidean'' TQFT. We show on examples that our
invariants are highly nontrivial.
\end{abstract}

\maketitle

\section{Introduction}

Let there be a Lie group~$G$ and its homogeneous space~$S$, and let the action of $G$ on~$S$ have an invariant depending on $k$ points, i.e., a function, called~$\Phi_2$ for further reasons, sending a $k$-tuple $(x_1,\dots,x_k)$ of points in~$S$ into an element of field $\mathbb F=\mathbb R \textrm{ or }\mathbb C$:
$$
\Phi_2\colon\; \underbrace{S\times\dots\times S}_{k\textrm{ times}}\to \mathbb F
$$
and such that its value does not change when an arbitrary element $g\in G$ acts on all~$x_1,\dots,x_k$:
$$
\Phi_2(gx_1,\dots,gx_k) = \Phi_2(x_1,\dots,x_k).
$$
Given a fixed $k$-tuple $(x_1^{(0)},\dots,x_k^{(0)})$ of points in~$S$, we consider also a mapping
$$\Phi_1\colon\; G\to \underbrace{S\times\dots\times S}_{k\textrm{ times}}
$$
given by formula
$$
\Phi_1\colon\; g\mapsto (gx_1^{(0)},\dots,gx_k^{(0)}).
$$
Then an obvious remark is that
\begin{equation}
\Phi_2\circ \Phi_1 = \const.
\label{Phi12}
\end{equation}
If we consider infinitesimal versions of $\Phi_1$ and~$\Phi_2$, i.e., tangent mappings $\varphi_1=d\Phi_1$ and~$\varphi_2=d\Phi_2$, the first of them taken at some arbitraty~$g$ and the second --- at the $k$-tuple $\Phi_1(g)$, then the consequence of~(\ref{Phi12}) is
\begin{equation}
\varphi_2\circ \varphi_1 =0.
\label{phi2phi1}
\end{equation}

It turns out that much more can be achieved if we have, in addition to $G$ and~$S$, a triangulated piecewise-linear manifold~$M$. Instead of just two mappings satisfying~(\ref{phi2phi1}), a meaningful \emph{chain complex} of vector spaces and their linear mappings $f_1,f_2,f_3,\dots$ can be constructed, at least for many specific $G$ and~$S$, such that $f_2\circ f_1=0$, $f_3\circ f_2=0$, and so on. The vector spaces consist of differentials of geometric values related to $G$, $S$ and the triangulation. Such chain complexes turn out to be acyclic in many cases, and the Reidemeister torsion for complexes of this kind can be used for constructing a wide range of manifold invariants.

Most of the work already done in this direction deals with the situation where $S=\mathbb R^3$ is a three-dimensional Euclidean space, $G=\mathrm E(3)$ is its group of motions and $M$ is a three-dimensional manifold. We will mention in this paper some of our works concerning this case; the latest achievement here was the construction of a finite-dimensional topological quantum field theory (TQFT)~\cite{K-tqft1,K-tqft2} in the spirit of M.~Atiyah's axioms~\cite{Atiyah}.

Other chain complexes for three-manifolds studied by us correspond to $S$ being an affine (real or complex) plane and $G$ --- the group of its
motions \emph{preserving the areas}~\cite{SL2-1,SL2-2}, and also to $S$ being a four-dimensional Euclidean space, $S=\mathbb R^4$, and $G=\mathrm
E(4)$ being its group of motions~\cite{man3geo4} (thus, both two- and four-dimensional homogeneous spaces~$S$ proved to be good for studying
\emph{three}-manifolds).

A chain complex for \emph{four}-manifolds has also been suggested where $S=\mathbb R^4$ and $G=\mathrm E(4)$~\cite{33,24,15}.

At this stage it is, however, too early to speak about a general recipe of how to construct a chain complex for an $m$-dimensional manifold, given a group~$G$ and its homogeneous space~$S$. Our current work consists rather in constructing and studying complexes for specific $G$ and~$S$. In the present paper, we investigate the case of $S=\mathbb C\cup\{\infty\}$ --- the compactified complex plane and $G=\mathrm{PSL}(2,\mathbb C)$ --- the group of its fractional-linear transformations. Our interest in this case was initially stimulated by the fact that it uses, as the reader will see below, some constructions known from hyperbolic geometry; it turned out later that there are also many new and beautiful features distinguishing this case from what was known earlier.

Below, in sections \ref{sec:cr} and~\ref{right} we construct what we call the basic complex --- a chain complex which is, in a sense, the simplest
possible one, and which is suitable for modifications used for various specific purposes. This construction goes in a somewhat unexpected way: we use
some geometric considerations in section~\ref{sec:cr} for constructing a half of the complex, and some rather different, at first sight,
considerations in section~\ref{right} for constructing its second half; the possibility to unite the two halves comes like a miracle. Then, in
section~\ref{sec:twisted} we construct a twisted version of the complex and prove its important property --- acyclicity. In section~\ref{sec:moves}
we produce manifold invariants using twisted complexes. In section~\ref{sec:examples} we provide some examples, together with one more --- relative
--- version of our complex in subsection~\ref{subsec:relative}. Finally, we discuss our results and further research in section~\ref{sec:discussion}.

\section{The left-hand half of the basic complex}
\label{sec:cr}

To begin, we consider a three-dimensional closed oriented manifold~$M$. We attach a complex number~$\zeta_i$ to every vertex~$i$ of its given
triangulation; $\zeta_i$ will be called the unperturbed, or initial, coordinate of vertex~$i$. These $\zeta$'s are parameters of our theory, of which
the final result will not depend. The only condition on $\zeta$'s is that they must lie in the general position with regard to all algebraic
constructions given below.

Now we define mappings~$F_1$, $F_2$ and~$F_3$. Mapping~$F_1$ sends an element of group~$\mathrm{PSL}(2,\mathbb C)$ represented by matrix $\begin{pmatrix} \alpha&\beta\\
\gamma&\delta
\end{pmatrix}$ into the column vector of height $N_0$ consisting of ``perturbed coordinates''
$$
z_i=\frac{\alpha\zeta_i+\beta}{\gamma\zeta_i+\delta}
$$
for all vertices~$i$; here $N_0$ is the number of vertices in the triangulation of~$M$.

The next mapping~$F_2$ sends a column vector of $N_0$ arbitrary values $z_i$ into a column vector of height~$N_3$, where $N_3$ is the number of
\emph{tetrahedra} in the triangulation. Each entry of this latter vector corresponds to a tetrahedron in the triangulation and is described as
follows. Let there be a tetrahedron~$0123$, whose orientation, given by this order of its vertices, corresponds to the given orientation of~$M$. The
entry of the mentioned vector, corresponding to tetrahedron~$0123$, consists of three complex values corresponding to its six \emph{unoriented} edges
and related as follows:
\begin{itemize}
\item the same value corresponds to any of two opposite edges: if $x$ corresponds to edge~$02$, it also corresponds to edge~$13$;
\item if $x$ corresponds to edges $02$ and~$13$, then the first of the values
\begin{equation}
1-\frac{1}{x}\,,\quad \frac{1}{1-x}
\label{2more}
\end{equation}
corresponds to any of the edges 03 and 12, while the second --- to the edges 01 and~23.
\end{itemize}
By definition, the~$x$ obtained by applying $F_2$ to given~$z$'s (where the actual tetrahedron vertices must be substituted instead of $0,1,2,3$) is
the cross-ratio
\begin{equation}
x=\frac{z_{01}z_{23}}{z_{03}z_{21}},
\label{cr}
\end{equation}
where
\begin{equation}
z_{ij}=z_i-z_j.
\label{zij}
\end{equation}
One can check that expressions (\ref{2more}) are in accordance with how the cross-ratio~(\ref{cr}) transforms under permutations of vertices.

Finally, mapping~$F_3$ sends a column vector of height~$N_3$ consisting of triples~$\bigl(x,\allowbreak\;1-1/x,\allowbreak\;1/(1-x)\bigr)$ into a
column vector of complex numbers~$\omega_{ij}$ of height~$N_1$, where $N_1$ is the number of edges in the triangulation, and $ij$ is a given edges
joining vertices $i$ and~$j$. Consider the \emph{star} of edge~$ij$; it consists of all tetrahedra having $ij$ as an edge. By definition, $F_3$
yields
\begin{equation}
\omega_{ij}=\prod x,
\label{oij}
\end{equation}
where all values $x$ in the product correspond to all tetrahedra in the star of~$ij$ and to the edge~$ij$ in each such tetrahedron. We call
$\omega_{ij}$ obtained according to formula \emph{deficit angle} around edge~$ij$.

Consider the following chain of spaces and mappings:
\begin{multline}
0 \longrightarrow \mathrm{PSL}(2,\mathbb C) \stackrel{\textstyle F_1}{\longrightarrow}
\left( \begin{array}{c} \textik{vertex }\\[-.5ex] \textik{coordinates}\\[-.5ex] z \end{array}\right)\\ \stackrel{\textstyle F_2}{\longrightarrow}
\left( \begin{array}{c} \textik{triples}\\ x,\,1-1/x,\,1/(1-x) \\ \textik{in tetrahedra}\end{array}\right) \stackrel{\textstyle F_3}{\longrightarrow}
\left( \begin{array}{c} \textik{deficit}\\[-.5ex] \textik{angles }\omega\\[-.5ex] \textik{around edges} \end{array}\right),
\label{crmacro}
\end{multline}
where the leftmost arrow, of course, just sends the zero into the unit of group $\mathrm{PSL}(2,\mathbb C)$.

\begin{theorem}
\label{th-l-glob}
The composition of any two successive arrows in~(\ref{crmacro}) is a constant mapping.
\end{theorem}

\begin{proof}
To show that $F_2\circ F_1=\const$, it is enough to say that the cross-ratio of four complex numbers is invariant under the action of the same element
of~$\mathrm{PSL}(2,\mathbb C)$ on all of them.

To show that $F_3\circ F_2=\const$, note that all terms in the product~(\ref{oij}) of values~(\ref{cr}) cancel out.
\end{proof}

We sometimes call the chain~(\ref{crmacro}) a ``macroscopic'' complex, in contrast to its differential, or ``microscopic'' version which we are going
to produce. Roughly speaking, it will consist of differentials of mappings $F_1$, $F_2$ and~$F_3$. This makes no difficulty when taking the
differential $f_1=dF_1\colon\; \mathfrak{psl}(2,\mathbb C)\to (dz)$, where $\mathfrak{psl}(2,\mathbb C)$ is the Lie algebra, and by $(dz)$ we denote
the vector space of column vectors of differentials of quantities~$z$. More formally, $(dz)$ is just a vector space over~$\mathbb C$ whose basis
consists of all the vertices of triangulation. To be exact, we represent $\mathfrak{psl}(2,\mathbb C)$ as a space of column vectors with components
$da,db,dc$, and define $f_1$ by the formula
\begin{equation}
dz_i = \begin{pmatrix}2\zeta_i & 1 & -\zeta_i^2 \end{pmatrix} \begin{pmatrix} da\\ db\\ dc \end{pmatrix}
\label{f1}
\end{equation}
for all vertices~$i$.

For the next mapping, we would like, however, to have if not one elegantly defined ``symmetric'' quantity instead of (\ref{cr}) and~(\ref{2more}),
then at least a simple differential not depending on the choice of an edge in the tetrahedron. For this, we propose
\begin{equation}
dy_{0123}=\frac{d\ln x}{\zeta_{02}\zeta_{31}}=\frac{d\ln (1-\frac{1}{x})}{\zeta_{03}\zeta_{12}}=\frac{d\ln \frac{1}{1-x}}{\zeta_{01}\zeta_{23}},
\label{dy}
\end{equation}
where $\zeta_{ij}=\zeta_i-\zeta_j$ similarly to~(\ref{zij}). Note, by the way, that such $dy$ does not even change even under {\em odd\/}
permutations of indices $0,1,2,3$. Thus, our mapping $f_2\colon\; (dz)\to (dy)$ is defined by differentiating formula~\ref{cr}; here $(dy)$ is the
space of column vectors whose coordinates are~$dy_{ijkl}$ for all tetrahedra~$ijkl$ in the triangulation or, more formally, a vector space
over~$\mathbb C$ whose basis consists of all the tetrahedra. The formula for~$f_2$ is:
\begin{equation}
dy_{0123}=\begin{pmatrix} -\frac{1}{\zeta_{01}\zeta_{02}\zeta_{03}} & -\frac{1}{\zeta_{01}\zeta_{12}\zeta_{13}} &
-\frac{1}{\zeta_{02}\zeta_{12}\zeta_{23}} & -\frac{1}{\zeta_{03}\zeta_{13}\zeta_{23}} \end{pmatrix} \begin{pmatrix} dz_0\\ dz_1\\ dz_2\\ dz_3
\end{pmatrix},
\label{f2}
\end{equation}
and similarly for other indices.

Finally, we introduce also logarithms of $\omega$'s in our definition of ``microscopic'' mapping~$f_3\colon\; (dy)\to (d\ln\omega)$, where
$(d\ln\omega)$ is again the obvious vector space, whose basis vectors are edges. We define $f_3$ by formula
\begin{equation}
d\ln\omega_{ij}=\zeta_{ij} \sum_{\mathrm{edges\ }kl} \zeta_{kl} \,dy_{ijkl},
\label{dlnomega}
\end{equation}
where ``edges $kl$'' are all edges belonging to the link of $ij$ and $ijkl$ gives the right orientation of every tetrahedron.

The following theorem is an obvious consequence of Theorem~\ref{th-l-glob}.

\begin{theorem}
\label{th-f}
The sequence of vector spaces and linear mappings
\begin{equation}
0\longrightarrow \mathfrak{psl}(2,\mathbb C) \stackrel{f_1}{\longrightarrow} (dz) \stackrel{f_2}{\longrightarrow} (dy)
\stackrel{f_3}{\longrightarrow} (d\ln\omega)
\label{crmicro}
\end{equation}
is a chain complex, i.e., $f_3\circ f_2=0$ and $f_2\circ f_1=0$.
\qed
\end{theorem}

Below in Section~\ref{right}, we will see how the complex~(\ref{crmicro}) can be continued to the right, and we will produce its ``right-hand half''.
Even after this, the resulting complex, literally taken, will be just the basis for building its different modifications useful for calculating
topological invariants and building a TQFT --- this is why we call (\ref{crmicro}) the left-hand half of the basic complex.

\section{The right-hand half of the basic complex and gluing the halves together}
\label{right}

Our next ``macroscopic'' sequence of spaces and (nonlinear) mappings is:
\begin{equation}
0 \longrightarrow \mathrm{SO}(3,\mathbb C) \stackrel{\textstyle G_1}{\longrightarrow}
\left( \begin{array}{c} \textik{isotropic}\\[-.5ex] \textik{vectors}\\[-.5ex] \textik{in vertices} \end{array}\right) \stackrel{\textstyle G_2}{\longrightarrow}
\left( \begin{array}{c} \textik{squared}\\[-.5ex] \textik{edge}\\[-.5ex] \textik{lengths} \end{array}\right) \stackrel{\textstyle G_3}{\longrightarrow}
\left( \begin{array}{c} \textik{discrepancies}\\ \Omega \\ \textik{in tetrahedra}\end{array}\right)
\label{crmacroconj}
\end{equation}

Here are the details. We consider a complex Euclidean space of column vectors of height~$3$ with the scalar product given by the matrix
\begin{equation}
\begin{pmatrix} 0&0&-1\\0&2&0\\-1&0&0 \end{pmatrix}.
\label{matrixEuclid}
\end{equation}
The group $\mathrm{SO}(3,\mathbb C)$ is realized as the group of matrices representing linear transformations of this space preserving the scalar product~(\ref{matrixEuclid}).

In each vertex~$i$ of the triangulation of manifold~$M$ now live \emph{two} complex parameters: $\zeta_i$ which is the same as in Section~\ref{sec:cr}, and a new parameter called~$\varkappa_i$. Of these, the following ``initial'', or unperturbed, \emph{isotropic vector} is made:
\begin{equation}
\vec e_i^{\textrm{ initial}} = \begin{pmatrix} \varkappa_i \zeta_i^2\\\varkappa_i \zeta_i\\ \varkappa_i \end{pmatrix}.
\label{isotropic-initial}
\end{equation}
The space called ``$\left( \begin{array}{c} \textik{isotropic}\\[-.5ex] \textik{vectors}\\[-.5ex] \textik{in vertices} \end{array}\right)$'' in~(\ref{crmacroconj}) consists of isotropic vectors~$\vec e_i$ in all vertices~$i$ of the form~(\ref{isotropic-initial}), but with all $\zeta_i$ and~$\varkappa_i$ replaced by arbitrary complex values $z_i$ and~$h_i$:
\begin{equation}
\vec e_i=\begin{pmatrix} h_i z_i^2\\h_i z_i\\h_i \end{pmatrix}
\label{ei}
\end{equation}
Our mapping~$G_1$ is defined as follows:
\begin{equation}
G_1\colon\quad T\in\mathrm{SO}(3,\mathbb C) \mapsto \{\textrm{vectors }\vec e_i = T \vec e_i^{\textrm{ initial}} \textrm{ for all }i\}.
\label{G1}
\end{equation}

The next space called ``$\left( \begin{array}{c} \textik{squared}\\[-.5ex] \textik{edge}\\[-.5ex] \textik{lengths} \end{array}\right)$'' in~(\ref{crmacroconj}) consists of complex numbers living on all \emph{edges}~$ij$. We assume that our isotropic vectors come out of the origin of coordinates. The map $G_2$ produces then, by definition, squared distances~$L_{ij}$ between the ends of $\vec e_i$ and~$\vec e_j$. Note the following relation with the scalar product:
\begin{equation}
L_{ij}=-2\vec e_i \vec e_j.
\label{Lij}
\end{equation}

Finally, our space ``$\left( \begin{array}{c} \textik{discrepancies}\\ \Omega \\ \textik{in tetrahedra}\end{array}\right)$'' consists of complex numbers~$\Omega_{ijkl}$ put in correspondence to all tetrahedra~$ijkl$. By definition, the $\Omega$'s produced by~$G_3$ from the given squared edge lengths are the following determinants:
\begin{equation}
\Omega_{ijkl}=\left| \begin{matrix}
0&L_{ij}&L_{ik}&L_{il}\\
L_{ji}&0&L_{jk}&L_{jl}\\
L_{ki}&L_{kj}&0&L_{kl}\\
L_{li}&L_{lj}&L_{lk}&0
\end{matrix} \right|,
\end{equation}
where of course $L_{ij}=L_{ji}$ and so on.

\begin{theorem}
\label{th-G-global}
The composition of any two successive arrows in~(\ref{crmacroconj}) is a constant mapping.
\end{theorem}

\begin{proof}
The relation $G_2\circ G_1 =\const$ holds simply because distances are invariant under the action of $\mathrm{SO}(3,\mathbb C)$.

The relation $G_3\circ G_2 =\mathrm{const}\; (=0)$ holds because $\Omega$ vanishes when the $L$'s are produced from \emph{three-dimensional} vectors according to~(\ref{Lij}).
\end{proof}

Now we pass on to ``microscopic'' values in full analogy with Section~\ref{sec:cr} and prove the following theorem --- analogue of Theorem~\ref{th-f} --- as an obvious consequence of Theorem~\ref{th-G-global}.

\begin{theorem}
\label{th-g}
The sequence of vector spaces and linear mappings
\begin{equation}
0\longrightarrow \mathfrak{so}(3,\mathbb C) \stackrel{\tilde g_1}{\longrightarrow} (dh)\oplus(dz) \stackrel{\tilde g_2}{\longrightarrow} (dL)
\stackrel{\tilde g_3}{\longrightarrow} (d\Omega),
\label{crmicroconj}
\end{equation}
where
$$
\tilde g_1=dG_1,\quad \tilde g_2=dG_2,\quad \tilde g_3=dG_3,
$$
is a chain complex, i.e., $\tilde g_3\circ \tilde g_2=0$ and $\tilde g_2\circ \tilde g_1=0$.
\qed
\end{theorem}

The tildes in Theorem~\ref{th-g} are because we actually want to slightly modify the complex~(\ref{crmicroconj}) by normalizing the bases in its vector spaces so as to make us able to join (\ref{crmicro}) and~(\ref{crmicroconj}) together in a single chain complex in the way described below. But first we must choose a basis in the Lie algebra~$\mathfrak{so}(3,\mathbb C)$. By definition, it will consist of the following three standard generators:
\begin{equation}
A=\begin{pmatrix}2&0&0\\0&0&0\\0&0&-2\end{pmatrix},\quad B=\begin{pmatrix}0&2&0\\0&0&1\\0&0&0\end{pmatrix},\quad
C=\begin{pmatrix}0&0&0\\1&0&0\\0&2&0\end{pmatrix}.
\end{equation}

Let $da^*,db^*,dc^*$ be infinitesimal numbers; here and below we mark with a star certain differentials, having in mind that we are going to ``transpose'' complex~(\ref{crmicroconj}) as described below, then the corresponding differentials without stars will enter in the transposed complex. We would like also to denote
\begin{equation}
d\alpha_i^*=\frac{dh_i}{2\varkappa_i},\quad d\beta_i^*=dz_i.
\label{alpha*beta*}
\end{equation}
If we calculate the change of $h_i$ and~$z_i$ under the action of matrix $da^*A+db^*B+dc^*C$ on vector $\vec e_i$~(\ref{ei}) and then substitute the initial values $h_i=\varkappa_i$ and $z_i=\zeta_i$ into the resulting Jacobian matrix, we get, taking also (\ref{alpha*beta*}) into account:
\begin{equation}
\begin{pmatrix}d\alpha_i^*\\d\beta_i^*\end{pmatrix}=\begin{pmatrix}-1&0&\zeta_i\\
2\zeta_i&1&-\zeta_i^2\end{pmatrix} \begin{pmatrix}da^*\\db^*\\dc^*\end{pmatrix}.
\label{g1}
\end{equation}
Formula~(\ref{g1}) gives our definition for linear mapping~$g_1$ --- the modified version of $\tilde g_1$ from~(\ref{crmicroconj}).

Next, we introduce normalized squared edge lengths in the following way:
\begin{equation*}
\varphi_{ij}^*=\frac{L_{ij}}{4 \varkappa_i \varkappa_j (\zeta_i-\zeta_j)^2}.
\end{equation*}
Thus, when $\varphi_{ij}^*$ is obtained according to $G_2$, it is
\begin{equation}
\varphi_{ij}^*=\frac{1}{2} \frac{h_i h_j (z_i-z_j)^2}{\varkappa_i \varkappa_j (\zeta_i-\zeta_j)^2}.
\label{eq:phiij}
\end{equation}
This yields
\begin{equation}
\frac{\partial \varphi_{ij}^*}{\partial \alpha_i^*}=1,\quad \frac{\partial \varphi_{ij}^*}{\partial \beta_i^*}=\frac{1}{\zeta_i-\zeta_j}.
\label{g2}
\end{equation}
By definition, formula~(\ref{g2}) gives matrix elements for linear mapping~$g_2$ --- the modified version of $\tilde g_2$ from~(\ref{crmicroconj}).

Finally, if $\Omega_{ijkl}$ is obtained according to~$G_3$ and we calculate the derivative $\partial \Omega_{ijkl}/\partial \varphi_{ij}^*$ at the
point where $L_{ij}=-2\vec e_i\vec e_j =2\varkappa_i\varkappa_i (\zeta_i-\zeta_j)^2$ and similarly for $L$'s with other indices, we get
\begin{equation*}
\frac{\partial \Omega_{ijkl}}{\partial \varphi_{ij}^*}=-128 (\zeta_i-\zeta_j)(\zeta_k-\zeta_l) \prod_{r<s} (\zeta_r-\zeta_s),
\end{equation*}
where in the product both $r$ and $s$ take values $i,j,k,l$, and ``$<$'' in ``$r<s$'' means just the alphabetic order. This suggests us to denote
\begin{equation*}
dy_{ijkl}^*=-\frac{d\Omega_{ijkl}}{128\prod_{r<s} (\zeta_r-\zeta_s)},
\end{equation*}
which yields
\begin{equation}
\frac{\partial y_{ijkl}^*}{\partial \varphi_{ij}^*} = \frac{1}{\zeta_{ij}\zeta_{kl}}.
\label{yphi}
\end{equation}
By definition, \eqref{yphi} gives matrix elements for linear mapping~$g_3$ --- the modified version of $\tilde g_3$ in~(\ref{crmicroconj}).

Hence, the modified version of~\eqref{crmicroconj} is
\begin{equation}
0\longrightarrow \mathfrak{so}(3,\mathbb C) \stackrel{g_1}{\longrightarrow} (d\alpha^*)\oplus(d\beta^*) \stackrel{g_2}{\longrightarrow} (d\varphi^*) \stackrel{g_3}{\longrightarrow} (dy^*),
\label{crmicroconj-modified}
\end{equation}

Comparing (\ref{yphi}) with (\ref{dlnomega}), we see that $f_3$ and~$g_3$ are related by matrix transposing:
\begin{equation}
g_3=f_3^{\mathrm T}.
\label{g3f3}
\end{equation}
This remarkable observation is the key for joining together our complexes \eqref{crmicro} and~\eqref{crmicroconj-modified}. Namely, here is our final basic complex:
\begin{multline}
0\longrightarrow \mathfrak{psl}(2,\mathbb C) \stackrel{f_1}{\longrightarrow} (dz) \stackrel{f_2}{\longrightarrow} (dy)
\stackrel{f_3}{\longrightarrow} (d\varphi)\\
\stackrel{f_4}{\longrightarrow} (d\alpha)\oplus(d\beta) \stackrel{f_5}{\longrightarrow} \mathfrak{psl}(2,\mathbb C)^* \longrightarrow 0.
\label{basic}
\end{multline}
By definition, in \eqref{basic}
\begin{equation*}
f_4=g_2^{\mathrm T},\quad f_5=g_1^{\mathrm T}.
\end{equation*}
As for the vector spaces, first, $(d\varphi)$ is just a new notation for the same~$(d\ln\omega)$ in~\eqref{crmicro}, which is justified by the fact that $\partial y_{ijkl}^* / \partial \varphi_{ij}^* = \partial (\ln\omega_{ij}) / \partial y_{ijkl}$, according to~\eqref{g3f3}. Next, $(d\alpha)$, $(d\beta)$ and~$\mathfrak{psl}(2,\mathbb C)^*$ can be considered just as convenient notations for some spaces of column vectors which are in an obvious sense dual to our spaces $(d\alpha^*)$, $(d\beta^*)$ and~$\mathfrak{so}(3,\mathbb C)$ respectively; in the latter case, we have taken into account the well-known isomorphism between Lie algebras.

\section{Twisted version of the complex}
\label{sec:twisted}

\subsection{Generalities on constructing the twisted version}

As we have already stated, we are going to use our basic complex~\eqref{basic} mainly as a starting point for various modifications. One important
modification is a complex twisted by a representation~$\rho$ of the fundamental group $\pi_1(M)$ of our considered manifold~$M$ into the
group~$\mathrm{PSL}(2,\mathbb C)\cong \mathrm{SO}(3,\mathbb C)$. This construction is similar to what we have done for the ``Euclidean''
case~\cite{twisted-euclidean1,twisted-euclidean2,M-thesis} and goes, in a few words, as follows:
\begin{enumerate}
\item \label{i1} we bring into consideration the {\em universal cover}~$\tilde M$ of~$M$, to whose vertices we now assign ``coordinate'' parameters in a way consistent with~$\rho$;
\item \label{i2} we add the parameters of the possible deformations of~$\rho$ to the second (nonzero) terms \emph{both from the left and from the right} in sequence~\eqref{basic};
\item \label{i3} we reduce the first terms both from the left and from the right --- the Lie algebra and its dual --- to the subalgebra commuting with the whole representation~$\rho$ and the dual space to that subalgebra.
\end{enumerate}
This all goes mostly in the same way as in the ``Euclidean'' case; below are some more details.

Item~(\ref{i1}): we assign coordinates $\varkappa$ and~$\zeta$ to the vertices of $\tilde M$ in a way consistent with representation~$\rho$: this
means that for any two vertices $i^{(1)},i^{(2)}\in \tilde M$ lying above the same vertex~$i\in M$, if $g\in\pi_1(M)$ is the element taking the first
of them into the second: $i^{(2)}=gi^{(1)}$, then
\begin{equation}
\begin{pmatrix} \varkappa_{i^{(2)}}\zeta_{i^{(2)}}^2\\ \varkappa_{i^{(2)}}\zeta_{i^{(2)}}\\ \varkappa_{i^{(2)}} \end{pmatrix} = \rho(g)
\begin{pmatrix} \varkappa_{i^{(1)}}\zeta_{i^{(1)}}^2\\ \varkappa_{i^{(1)}}\zeta_{i^{(1)}}\\ \varkappa_{i^{(1)}} \end{pmatrix}.
\label{ini-hz}
\end{equation}
Here $\rho(g)$ is understood as the element of $\mathrm{SO}(3,\mathbb C)$ in the sense of Section~\ref{right}; the corresponding element of the isomorphic group~$\mathrm{PSL}(2,\mathbb C)$ is the fractional-linear transformation of $\zeta$'s determined by~\eqref{ini-hz}. Otherwise, the coordinates are arbitrary with the only condition of general position with respect to all our algebraic constructions.

When we (slightly) deform these ``initial'' coordinates, their respective values must still obey the same restriction as~(\ref{ini-hz}), i.e.,
\begin{equation}
\begin{pmatrix} h_{i^{(2)}}z_{i^{(2)}}^2\\ h_{i^{(2)}}z_{i^{(2)}}\\ h_{i^{(2)}} \end{pmatrix} = \rho(g)
\begin{pmatrix} h_{i^{(1)}}z_{i^{(1)}}^2\\ h_{i^{(1)}}z_{i^{(1)}}\\ h_{i^{(1)}} \end{pmatrix}.
\label{hz}
\end{equation}
This means that, in the ``microscopic'' complex~\eqref{basic}, the differentials $dh$ and~$dz$ for different copies of the same vertex~$i$ can be obtained one from another by differentiating formula~(\ref{hz}). So, the second nonzero terms from both sides in~\eqref{basic} are now modified as follows: they contain the differentials $dh$ and~$dz$ (the left one just~$dz$, of course) for \emph{one} ``main'' copy (lying in~$\tilde M$ and arbitrarily chosen) of each vertex~$i\in M$. According to the general Definition~\ref{dfn:family} below in Section~\ref{sec:moves}, these ``main'' copies form a \emph{fundamental family} of vertices in~$\tilde M$.

Item~(\ref{i2}): we assume that $\rho$ is a {\em regular\/} point in the space of all representations $\pi_1(M)\to \mathrm{SO}(3,\mathbb C)$ in the sense that its neighborhood --- all near representations --- can be parameterized (smoothly enough) by some number of parameters. To be exact, by ``representation'' we understand here a \emph{class of equivalent} representations.
We assume also that we are able to choose one specific representative in each class of equivalent representations, and these representatives are also smoothly parameterized by the same parameters. Such a specific representative is needed to calculate the deformed coordinates of all copies of vertices in~$\tilde M$ from their ``main'' copies. This is done, of course, according to the same formula~(\ref{hz}) but now with the deformed~$\rho$.

One more modification to the second nonzero term from the left in~\eqref{basic} is adding there the vector space~$(dg)$ of all the infinitesimal
representation deformation parameters; the corresponding change for the second nonzero term from the right is adding the dual space~$(dg)^*$.

Item~(\ref{i3}): we denote the subalgebra of $\mathfrak{psl}(2,\mathbb C)$ commuting with~$\rho$ as $\mathfrak{psl}(2,\mathbb C)_{\rho}$; the $\rho$ here is of course undeformed.

Here is how we write the resulting twisted complex:
\begin{multline}
0\longrightarrow \mathfrak{psl}(2,\mathbb C)_{\rho} \stackrel{f_1}{\longrightarrow} (dz)\oplus (dg) \stackrel{f_2}{\longrightarrow} (dy) \stackrel{f_3}{\longrightarrow} (d\varphi)\\
\stackrel{f_4}{\longrightarrow} (d\alpha)\oplus(d\beta)\oplus (dg)^* \stackrel{f_5}{\longrightarrow} \mathfrak{psl}(2,\mathbb C)_{\rho}^* \longrightarrow 0.
\label{twisted}
\end{multline}

\subsection{Acyclicity of the twisted complex}
\label{sec:cr-acyclicity}

\begin{theorem}
\label{th-acyclicity}
Complex~\eqref{twisted} is acyclic for any closed oriented triangulated three-dimensional manifold~$M$, i.e., sequence~\eqref{twisted} is \emph{exact}: the image of any mapping coincides with the kernel of the next mapping.
\end{theorem}

\begin{remark}
Below in Subsection~\ref{acyclic-generalities}, we give also the definition of acyclicity from a homological viewpoint, see Definition~\ref{dfn:acycl} .
\end{remark}

The rest of this subsection contains the proof of Theorem~\ref{th-acyclicity}. It is divided into two parts: acyclicity in the first three terms (from the left) in~\eqref{twisted}, i.e., before the arrow~$f_3$, and acyclicity in the second three terms.

\subsubsection{Acyclicity in the first three terms}
\label{subsubsec:left}

This proof goes in a direct analogy with the ``Euclidean'' case, see~\cite{M-thesis} for the most detailed exposition.

Namely, first we consider the twisted analogue of the macroscopic complex~\eqref{crmacro}:
\begin{multline}
0 \longrightarrow \mathrm{PSL}(2,\mathbb C)_{\rho} \stackrel{\textstyle F_1}{\longrightarrow}
\left( \begin{array}{c} \textik{vertex coordinates } z\\[-.5ex] \textik{and deformation}\\[-.5ex] \textik{parameters }g \end{array}\right)\\ \stackrel{\textstyle F_2}{\longrightarrow}
\left( \begin{array}{c} \textik{triples}\\ x,\,1-1/x,\,1/(1-x) \\ \textik{in tetrahedra}\end{array}\right) \stackrel{\textstyle F_3}{\longrightarrow}
\left( \begin{array}{c} \textik{deficit}\\[-.5ex] \textik{angles }\omega\\[-.5ex] \textik{around edges} \end{array}\right),
\label{crmacro-twisted}
\end{multline}
Here $\mathrm{PSL}(2,\mathbb C)_{\rho}$ means the subgroup commuting with~$\rho$; parameters~$g$ smoothly parameterize the equivalence classes of representations in some neighborhood of~$\rho$, we assume that $\rho$ itself corresponds to $g=0$. Vertex coordinates~$z$ belong to a fundamental family of vertices. We want to prove that
\begin{enumerate}
\item \label{j1} the pre-image $F_1^{-1}(\zeta,0)$ of undeformed $z$'s and~$g$'s is exactly the unit of~$\mathrm{PSL}(2,\mathbb C)_{\rho}$,
\item \label{j2} if all the cross-ratios~$x$ obtained due to~$F_2$ remain constant, than the representation~$\rho$ stays in its equivalence class --- deformation parameters~$g=0$, and vertex coordinates~$z$ change from their initial values~$\zeta$ in such way that they can be obtained from some $\mathrm{PSL}(2,\mathbb C)_{\rho}$ element according to the mapping~$F_1$,
\item \label{j3} if values~$x$ in tetrahedra are such that all ``deficit angles''~$\omega$, obtained due to~$F_3$, are unities, then these~$x$'s can be obtained from some vertex coordinates~$z$ and deformation parameters~$g$ according to~$F_2$.
\end{enumerate}

Item (\ref{j1}) is clear. To prove (\ref{j2}), we note that the constancy of cross-ratios guarantees that the $z$'s \emph{in the whole universal cover}~$\tilde M$ can be obtained from the~$\zeta$'s by a single transformation $\gamma\in \mathrm{PSL}(2,\mathbb C)$; it easily follows from here that $\rho$ stays in its equivalence class, which is only possible in our situation if $\rho$ just remains the same and $\gamma$ commutes with~$\rho$.

So, the key issue here is to prove (\ref{j3}). We begin with choosing coordinates~$z$ for one arbitrarily chosen ``initial'' tetrahedron. Coordinates are arbitrary except that they must have the required cross-ratio~$x$. We say then that we have associated a \emph{coordinate system} with this initial tetrahedron. Then, we extend this coordinate system to tetrahedra having a common 2-face with the initial tetrahedron, which goes in a unique way given the $x$'s and these adjacent tetrahedra, and continue this procees to the next adjacent tetrahedra and so on. The vanishing deficit angles guarantee that the thus obtained coordinate system in any tetrahedron does not depend on the specific way joining it with the initial tetrahedron (recall also that we are in the simply connected universal cover~$\tilde M$). From the coordinates of vertices in~$\tilde M$ we can also extract the (deformed) representation $\pi_1(M)\to \mathrm{PSL}(2,\mathbb C)$.

Second comes the transition from macroscopic to microscopic situation, i.e., to the left-hand half of complex~\eqref{twisted}. This uses the fact
that the initial coordinates~$\zeta$ are in general position and goes quite similarly to such transition in the ``Euclidean'' situation
of~\cite{M-thesis}.

\subsubsection{Acyclicity in the second three terms}

Again, first comes the macroscopic part --- we consider the twisted version of sequence~\eqref{crmacroconj}:
\begin{multline}
0 \longrightarrow \mathrm{SO}(3,\mathbb C)_{\rho} \stackrel{\textstyle G_1}{\longrightarrow}
\left( \begin{array}{c} \textik{isotropic vectors in vertices}\\[-.5ex] \textik{and deformation}\\[-.5ex] \textik{parameters }g \end{array}\right)\\ \stackrel{\textstyle G_2}{\longrightarrow}
\left( \begin{array}{c} \textik{squared}\\[-.5ex] \textik{edge}\\[-.5ex] \textik{lengths} \end{array}\right) \stackrel{\textstyle G_3}{\longrightarrow}
\left( \begin{array}{c} \textik{discrepancies}\\ \Omega \\ \textik{in tetrahedra}\end{array}\right)
\label{crmacroconj-twisted}
\end{multline}

Again, the key point belongs, like in Subsubsection~\ref{subsubsec:left}, to the term ``(squared edge lengths)'': we must prove that, if edge lengths
are such that the discrepancies in all tetrahedra vanish, then these edge lengths can be obtained as distances between the ends of isotropic vectors
starting at the origin of coordinates, with a duly deformed representation~$\rho$. We begin again with choosing one ``initial'' tetrahedron and
assigning coordinates $\varkappa$ and~$\zeta$ to its vertices, which can be done due to the following lemma whose proof is a simple exercise in
linear algebra.

\begin{lemma}
\label{lemma:initial}
If all distances between the vertices of a tetrahedron are given, such that the discrepancy in this tetrahedron vanishes, then the vertices of this tetrahedron can be placed at the ends of isotropic vectors, and this is done uniquely up to an orthogonal rotation.
\qed
\end{lemma}

Then we proceed like in~\cite{M-thesis} and our Subsubsection~\ref{subsubsec:left}: we extend our coordinates to neighboring tetrahedra, i.e., having a common two-face with a tetrahedron whose vertices have already been assigned coordinates. Let such a new tetrahedron be~$1234$, with the mentioned common two-face~$123$. This means that we assign coordinates $\varkappa_4$ and~$\zeta_4$ to just one new vertex~$4$, for which three distances --- lengths of edges $14$, $24$ and~$34$ --- are given, with the condition of zero discrepancy.

\begin{lemma}
\label{lemma:extending}
Suppose that the coordinated $\varkappa$ and~$\zeta$ are given for vertices $1$, $2$ and~$3$, and the lengths of edges $14$, $24$ and~$34$ are given, with the condition of zero discrepancy in tetrahedron~$1234$. Then, the coordinates of vertex~$4$ are determined uniquely, if everything happens in a general position.
\end{lemma}

\begin{proof}
The ends of isotropic vectors form a cone $w_2^2=w_1w_3$ in the three-dimen\-sional complex Euclidean space of points $(w_1,w_2,w_3)$. If a point in
this cone is given, then a ``circle'' centered at this point, i.e., the set of points situated at some fixed distance from it, is a \emph{parabola}.
Two such generic parabolas intersect at two points, much like two circles in a two-dimensional sphere. So, if points $1$ and~$2$ are given together
with the lengths of edges $14$ and~$24$, then there are two possibilities for placing point~$4$. Only one of these possibilities is selected if we
have also a point~$3$ and the length of~$34$, with the compatibility condition of zero discrepancy.
\end{proof}

What remains is to show that the resulting coordinates do not depend on the way joining the initial and the final tetrahedron. Like
in~\cite{M-thesis} and Subsubsection~\ref{subsubsec:left}, it is enough to prove this fact just for the star of some edge in the triangulation, and it can be formulated as the following Lemma.

\begin{lemma}
If all the lengths for all edges in the star of some edge are given, with the conditions that
\begin{enumerate}
\item the discrepancy in each tetrahedron is zero,
\item there are some ``initial'' coordinates $\varkappa$ and~$\zeta$ for all vertices in the star and thus some initial values for the edge lengths, and the actual (``deformed'') edge lengths are close to the initial ones,
\end{enumerate}
then the vertices of the whole star can be assigned coordinates $h$ and~$z$ compatible with these lengths. This can be done up to an orthogonal rotation.
\end{lemma}

\begin{proof}
Consider the star of some edge~$12$. We first fix the coordinates $h$ and~$z$ for vertices $1$ and~$2$ at some points in the cone such that $h_k\approx \varkappa_k$,\quad $z_k\approx \zeta_k$,\quad $k=1,2$, and such that the distance between $1$ and~$2$ equals the required length of~$12$. Then, to determine the coordinates of any vertex~$i$ in the \emph{link} of~$12$, it is enough to know the required lengths of edges $i1$ and~$i2$ together with the fact that $i$ must lie close to its initial position (this latter condition throws away the unwanted possibility for the position of~$i$, see the proof of Lemma~\ref{lemma:extending} above). The right lengths of edges in the link of~$12$ are ensured automatically.
\end{proof}

As in Subsubsection~\ref{subsubsec:left}, the transition to microscopic case using general position argument completes the proof of acyclicity for
the right-hand half of~\eqref{basic} and thus for the whole complex~\eqref{twisted}.

\section{Torsion, Pachner moves, and a manifold invariant}
\label{sec:moves}

\subsection{Generalities on acyclic complexes and their torsions}
\label{acyclic-generalities}

The key value which we want to extract from complex~\eqref{basic} and similar algebraic complexes is its (Reidemeister) \emph{torsion}. To introduce this important notion properly, and for the reader's convenience, we remind here briefly basic definitions from the theory of algebraic complexes, including those already used in this paper. More detailed exposition can be found in monograph~\cite{Tur01}.

Let $C_0$, $C_1$, \ldots , $C_n$ be finite-dimensional $\mathbb{C}$-vector spaces. We suppose that each $C_i$ is based, that is, a distinguished
basis in it is indicated. Then, a linear mapping $f_i \colon\; C_{i + 1} \to C_i$ can be identified with a matrix.

\begin{dfn}
\label{dfn:cmplx}
The sequence of vector spaces and linear mappings
\begin{equation}\label{eq:cmplx}
C = (0 \xrightarrow{} C_n \xrightarrow{f_{n-1}} C_{n-1} \xrightarrow{} \ldots \xrightarrow{} C_1 \xrightarrow{f_0} C_0 \xrightarrow{} 0)
\end{equation}
is called a \textit{complex} if $\Ima f_i \subset \Ker f_{i-1}$ for all $i = 1, \ldots , n - 1$. This condition is equivalent to $f_{i-1} \, f_i = 0$ for all $i$.
\end{dfn}

\begin{dfn}
\label{dfn:homol}
The space $H_i(C) = \Ker f_{i-1} / \Ima f_i$ is called the \textit{$i$th homology} of the complex $C$.
\end{dfn}

\begin{dfn}\label{dfn:acycl}
The complex $C$ is said to be \textit{acyclic} if $H_i(C) = 0$ for all $i$. This condition is equivalent to $\rank f_{i-1} = \dim C_i - \rank f_i$ for all $i$.
\end{dfn}

Suppose the sequence~\eqref{eq:cmplx} is an acyclic complex. Let $\mathcal{C}_i$ be an ordered set of basis vectors in $C_i$ and $\mathcal{B}_i
\subset \mathcal{C}_i$ be a subset of basis vectors belonging to the space~$\Ima f_i$.

Denote by ${}_{\mathcal{B}_i} f_i$ a nondegenerate transition matrix from the basis in space $C_{i+1}/\Ima f_{i+1}$ to the basis in space $\Ima f_i$.
By acyclicity, such a matrix does exist. Hence, ${}_{\mathcal{B}_i} f_i$ is a principal minor of the matrix $f_i$ obtained by striking out the rows
corresponding to vectors of $\mathcal{B}_{i+1}$ and the columns corresponding to vectors of $\mathcal{C}_i \setminus \mathcal{B}_i$.

\begin{dfn}\label{dfn:tors}
The quantity
\begin{equation}
\label{eq:torsion} \tau(C) = \prod\limits_{i = 0}^{n-1} (\det {}_{\mathcal{B}_i} f_i)^{(- 1)^i}
\end{equation}
is called the \textit{torsion} of acyclic complex $C$.
\end{dfn}

\begin{remark}
The torsion $\tau(C)$ defined above is the inverse of the torsion defined in~\cite{Tur01}.
\end{remark}

\begin{theorem}[\cite{Tur01}]
Up to a sign, $\tau(C)$ does not depend on the choice of subsets~$\mathcal{B}_i$.
\end{theorem}

\begin{remark}\label{rmk:changing}
The torsion $\tau(C)$ does depend on the distinguished basis of $C_i$. If one performs change-of-basis transformation in every space $C_i$ with
nondegenerate matrix $A_i$, then the torsion $\tau(C)$ is multiplied by
\[
\prod\limits_{i = 0}^n (\det A_i)^{(-1)^i}.
\]
\end{remark}

Let us define a nondegenerate $\tau$-chain following V.~Turaev~\cite{Tur01}.
\begin{dfn}
\label{dfn:tauchain}
Let $\alpha_i$ be certain collection of basis vectors in the space $C_i$ of~\eqref{eq:cmplx}. Let $S_i$ be a submatrix of~$f_i$ generated by such elements $a^i_{jk}$ that $j$ corresponds to some element from $\alpha_{i+1}$ and $k$ corresponds to some element from $\overline{\alpha}_i$. A collection of sets $(\alpha_0, \alpha_1, \ldots, \alpha_n)$ is called a \textit{nondegenerate $\tau$-chain} if the matrices $S_i$ are square and nondegenerate for all~$i$.
\end{dfn}
\begin{lemma}[\cite{Tur01}]
\label{lem:tauchain}
Complex~\eqref{eq:cmplx} is acyclic if and only if it has a nondegenerate $\tau$-chain.
\end{lemma}

\subsection{Invariant}

There are different invariants obtained from various modification of the basic complex~\eqref{basic}, including invariants of manifolds with
boundary, knots and links. Here we are going to consider a twisted complex~\eqref{twisted} for a closed oriented manifold~$M$, having in mind that
all reasonings of this subsection are easily modified for other situations. Note that our linear mappings~$f_1,\dots,f_5$ in~\eqref{twisted} are
numbered in a different way compared to~\eqref{eq:cmplx}; this does not bring about any serious changes.

Denote by $\mathcal{C}_i$ ($i = 0,\ldots, 5$) arbitrary ordered sets of basis vectors in all the spaces of~\eqref{twisted} starting from
$\mathfrak{psl}(2, \mathbb{C})_\rho$. Let $\mathcal{B}_i \subset \mathcal{C}_i$ be a subset of basis vectors belonging to $\Ima f_i$. Denote by
${}_{\mathcal{B}_i} f_i$ such a principal minor of the matrix $f_i$ that its rows correspond to the vectors from $\mathcal{C}_{i-1} \setminus
\mathcal{B}_{i-1}$ and its columns correspond to the vectors from $\mathcal{B}_i$. Due to acyclicity of~\eqref{twisted}, ${}_{\mathcal{B}_i} f_i$ is
really exists. Suppose that $\mathcal{C}_4 \setminus \mathcal{B}_4 = \mathcal{B}_1$. Then, according to Definition~\ref{dfn:tors}, the torsion of
complex~\eqref{twisted} looks like
\begin{equation}\label{eq:tors}
\tau = \frac{(\det {}_{\mathcal{B}_1} f_1)^2 \, \det {}_{\mathcal{B}_3} f_3}{\det {}_{\mathcal{B}_2} f_2 \, \det {}_{\mathcal{B}_4} f_4}.
\end{equation}

Recall that $\tilde{M}$ denotes the universal cover of our triangulated 3-manifold~$M$.

\begin{dfn}\label{dfn:family}
A \emph{fundamental family of simplices} in $\tilde{M}$ is such a family $\mathcal{F}$ of simplices of $\tilde{M}$ that over each simplex of $M$ lies
exactly one simplex of this family.
\end{dfn}

\begin{theorem}\label{thm:inv}
The quantity
\begin{equation}\label{eq:inv}
I_\rho(M; \mathrm{PSL}(2, \mathbb{C})) =  \frac{\tau}{\prod \zeta^2_{ij}}
\end{equation}
is a topological invariant of manifold~$M$. Here $\zeta_{ij} = \zeta_i - \zeta_j$ for the edge $ij$ and the product is taken over all edges from the
fundamental family~$\mathcal{F}$.
\end{theorem}

\begin{remark}
\label{remark_sign} As it is known (see monograph~\cite{Tur01}), usually a torsion is defined up to a sign, so that special measures must be taken
for its ``sign-refining''. This sign is changed when we change the order of basis vectors in any of the vector spaces. In the present paper, we
assume that the value~\eqref{eq:inv} and other similar values below are taken up to a sign.
\end{remark}

\begin{proof}
We are going to us show that $I_\rho(M; \mathrm{PSL(2, \mathbb{C})})$ is invariant under the Pachner moves $2 \to 3$ and~$1 \to 4$.

Recall that a move $2 \to 3$ replaces two adjacent concordantly oriented tetrahedra $1234$ and~$5123$ with three tetrahedra $1254$, $2354$ and $3154$
by adding a new edge~$45$ into the triangulation. Let us denote by $\tilde{\mathcal{B}}_i$ the set of basis vectors from $\Ima f_i$ after doing the
move~$2 \to 3$.

We set $\tilde{\mathcal{B}}_2 = (\mathcal{B}_2 \setminus \{dy_{1234}, dy_{5123}\}) \cup \{dy_{2354}, dy_{3154}\}$. Then, applying the row expansion
of the determinants $\det {}_{\tilde{\mathcal{B}}_2} \tilde{f}_2$ and $\det {}_{\mathcal{B}_2} f_2$ and comparing the corresponding multipliers, one
can see that
\[
\det {}_{\tilde{\mathcal{B}}_2} \tilde{f}_2 = \det {}_{\mathcal{B}_2} f_2 \cdot \frac{\zeta_{12}}{\zeta_{54}}.
\]

Further, using the same argumentation as in the Euclidean case (\cite[Lemma 1.15]{M-thesis}), we get
\[
\det {}_{\tilde{\mathcal{B}}_3} \tilde{f}_3 = \det {}_{\mathcal{B}_3} f_3 \cdot \pa{\varphi_{45}}{y_{1254}} = \det {}_{\mathcal{B}_3} f_3 \cdot
\zeta_{12}\,\zeta_{54},
\]
where $\tilde{\mathcal{B}}_3 = \mathcal{B}_3 \cup \{d\varphi_{45}\}$ and formula~\eqref{dlnomega} is taken into account. So,
\begin{equation}\label{eq:tt}
\frac{\tilde{\tau}}{\tau} = \zeta_{54}^2,
\end{equation}
where $\tau$ is the torsion of~\eqref{twisted} before doing the move $2 \to 3$ and $\tilde{\tau}$ is the one after. It follows that
quantity~\eqref{eq:inv} does not change under the move $2 \to 3$. Besides, by the lemma on nondegenerate $\tau$-chain (Lemma~\ref{lem:tauchain}), it
follows from~\eqref{eq:tt} that the complex~\eqref{twisted} remains acyclic.

Consider now a move $1 \to 4$. A new vertex $5$ is added into a tetrahedron $1234$, which is replaced so with four concordantly oriented tetrahedra $1235$, $1254$, $2354$
and~$3154$.

We set $\tilde{\mathcal{B}}_2 = (\mathcal{B}_2 \setminus \{dy_{1234}\}) \cup \{dy_{2354}, dy_{3154}\}$. Then, again applying the row expansion, it is
easy to show that
\[
\det {}_{\tilde{\mathcal{B}}_2} \tilde{f}_2 = \det {}_{\mathcal{B}_2} f_2 \cdot \frac{\zeta_{12}}{\zeta_{15}\,\zeta_{25}\,\zeta_{35}\,\zeta_{54}}.
\]

Let us choose $\tilde{\mathcal{B}}_3 = \mathcal{B}_3 \cup \{d\varphi_{CE}, d\varphi_{DE}\}$. Then, using formula~\eqref{dlnomega}, we have
\[
\det {}_{\tilde{\mathcal{B}}_3} \tilde{f}_3 = \det {}_{\mathcal{B}_3} f_3 \cdot
\begin{vmatrix}
\pa{\varphi_{35}}{y_{1235}} & \pa{\varphi_{35}}{y_{1254}} \\[3mm]
\pa{\varphi_{45}}{y_{1235}} & \pa{\varphi_{45}}{y_{1254}}
\end{vmatrix} = \det {}_{\mathcal{B}_3} f_3 \cdot \zeta_{12}^2\,\zeta_{35}\,\zeta_{54}.
\]

Similarly, if $\tilde{\mathcal{B}}_4 = \mathcal{B}_4 \cup \{d\alpha_5, d\beta_5\}$, then
\[
\det {}_{\tilde{\mathcal{B}}_4} \tilde{f}_4 = \det {}_{\mathcal{B}_4} f_4 \cdot
\begin{vmatrix}
\pa{\alpha_5}{\varphi_{15}} & \pa{\alpha_5}{\varphi_{25}} \\[3mm]
\pa{\beta_5}{\varphi_{15}} & \pa{\beta_5}{\varphi_{25}}
\end{vmatrix} = \det {}_{\mathcal{B}_4} f_4 \cdot
\begin{vmatrix}
1 & 1 \\[3mm]
\frac{1}{\zeta_{15}} & \frac{1}{\zeta_{25}}
\end{vmatrix} = \det
{}_{\mathcal{B}_4} f_4 \cdot \frac{\zeta_{12}}{\zeta_{15}\,\zeta_{25}},
\]
where we have used the fact that the partial derivatives in the minor are the same as the corresponding elements of the transposed matrix
$g_2=f_4^{\mathrm T}$, see formulas~\eqref{g2}.

Combining the above results, we obtain
\begin{equation}\label{eq:tt2}
\frac{\tilde{\tau}}{\tau} = (\zeta_{15}\,\zeta_{25}\,\zeta_{35}\,\zeta_{54})^2.
\end{equation}
It follows that~\eqref{eq:inv} is invariant under the move $1 \to 4$. Again, by the lemma on nondegenerate $\tau$-chain, the complex~\eqref{twisted}
remains acyclic.

Recall that, by Pachner's theorem~\cite{Pach}, invariance under the Pachner moves means topological invariance.

We must also show that $I_\rho(M; \mathrm{PSL(2, \mathbb{C})})$ does not depend on any detail of its construction. In particular, it is invariant
under the choice of fundamental family~$\mathcal{F}$ and initial values~$\zeta_i$. One can easily check this by analogy with the Euclidean case
(see~\cite[Theorem 1.13]{M-thesis} for details). Theorem~\ref{thm:inv} is proven.

\end{proof}

\section{Examples}
\label{sec:examples}

\subsection{$\boldsymbol{S^2\times S^1}$}
\label{S2xS1}

We consider a triangulation of~$S^2 \times S^1$ consisting of two copies of the same triangular prism with bases $123$ and~$1'2'3'$. Each of these
prisms is divided in three tetrahedra $11'23$, $1'22'3'$ and~$1'233'$ by adding three edges $1'2$, $1'3$ and~$23'$.

Let $\eta$ denote a generator of $\pi_1(S^2 \times S^1) \cong \mathbb{Z}$. Below we consider two possibilities for a representation $\rho \colon
\;\mathbb{Z} \to \mathrm{PSL}(2, \mathbb{C})$.

\subsubsection{Non-parabolic representation}

Here $\rho$ is given by
\begin{equation}
\rho \colon \; \eta \mapsto \begin{pmatrix} \lambda & 0 \\ 0 & \lambda^{-1} \end{pmatrix},
\end{equation}
where $\lambda \neq 0, 1$. Then, $(\rho \eta)z = \lambda^2 z$. Since $\varphi_{ij}^*$ from~\eqref{eq:phiij} must be invariant, it follows that $(\rho
\eta)h = \lambda^{-2}h$.

The algebra $\mathfrak{psl}(2, \mathbb{C})_\rho$ in~\eqref{twisted} is 1-dimensional, its basis consists of~$da$, see formula~\eqref{f1}. Therefore,
$\mathcal{C}_0 = \mathcal{B}_0 = \{da\}$.

The space $(dg)$ is also 1-dimensional. For reasons of symmetry between the two possible generators of $\mathbb{Z}$, we assume that $(dg)$ is
generated by~$d\lambda/\lambda$.

The space $(dz)$ is three-dimensional and generated by $dz$'s of vertices $1$, $2$ and~$3$. Let us choose
\[
\mathcal{B}_1 = \{dz_{1}\}.
\]
Then, obviously,
\[
\det {}_{\mathcal{B}_1} f_1 = 2\zeta_1.
\]

Choosing $\mathcal{B}_3 = \{d\varphi_{11'}, d\varphi_{22'}, d\varphi_{33'}\}$, we get
\begin{multline}
\det {}_{\mathcal{B}_3} f_3 =
\begin{vmatrix}
\pa{\varphi_{11'}}{y_{1'123}} & \pa{\varphi_{11'}}{y_{1'2'23'}} & \pa{\varphi_{11'}}{y_{1'23'3}} \\[3mm]
\pa{\varphi_{22'}}{y_{1'123}} & \pa{\varphi_{22'}}{y_{1'2'23'}} & \pa{\varphi_{22'}}{y_{1'23'3}} \\[3mm]
\pa{\varphi_{33'}}{y_{1'123}} & \pa{\varphi_{33'}}{y_{1'2'23'}} & \pa{\varphi_{33'}}{y_{1'23'3}}
\end{vmatrix}  \\[.5ex]
= \begin{vmatrix}
\zeta_{1'1} \zeta_{23} & 0 & 0 \\
0 & \zeta_{2'2} \zeta_{1'3'} & 0 \\
0 & 0 & \zeta_{3'3} \zeta_{1'2}
\end{vmatrix} = \lambda^2 \, \zeta_{1'1} \, \zeta_{2'2} \, \zeta_{3'3} \, \zeta_{1'2} \, \zeta_{13} \, \zeta_{23}.
\end{multline}

Similarly, choosing $\mathcal{B}_2 = \{dy_{11'23}, dy_{1'22'3'}, dy_{1'233'}\}$ and $\mathcal{B}_4 = \{d\alpha_1, d\alpha_2, \allowbreak d\alpha_3,
\allowbreak d\beta_2, \allowbreak d\beta_3, (d\lambda/\lambda)^*\}$, we have calculated
\[
\det {}_{\mathcal{B}_2} f_2 = - \frac{2(\lambda - 1)^2 \, (\lambda + 1)^2 \, \zeta_1}{\lambda^2 \, \zeta_{1'1} \, \zeta_{2'2} \, \zeta_{3'3} \,
\zeta_{12} \, \zeta_{13} \, \zeta_{23} \, \zeta_{1'2} \, \zeta_{1'3} \, \zeta_{23'}}
\]
and
\[
\det {}_{\mathcal{B}_4} f_4 = \frac{2(\lambda - 1)^2 \, (\lambda + 1)^2 \, \zeta_1}{\zeta_{12} \, \zeta_{23'} \, \zeta_{1'3}}.
\]

Substituting all the found minors into~\eqref{eq:tors}, we obtain the torsion of~\eqref{twisted} and then, by formula~\eqref{eq:inv}, find
\begin{equation}\label{eq:IS2S1nonparab}
I_\rho(S^2 \times S^1; \mathrm{PSL(2, \mathbb{C})}) = - \left(\lambda - \frac{1}{\lambda}\right)^{-4}.
\end{equation}

\subsubsection{Parabolic representation}

Here $\rho \colon\; \mathbb{Z} \to \mathrm{PSL}(2, \mathbb{C})$ is defined by
\begin{equation}\label{eq:parab}
\eta \mapsto \begin{pmatrix} 1 & 1 \\ 0 & 1 \end{pmatrix}.
\end{equation}
Then, $(\rho\eta) z = z+1$. Again, since $\varphi_{ij}^*$ from~\eqref{eq:phiij} must be invariant, we have $(\rho\eta)h = h$.

The algebra $\mathfrak{psl}(2, \mathbb{C})_\rho$ in~\eqref{twisted} is again 1-dimensional, but now its basis consists of $db$, see again
formula~\eqref{f1}. Therefore, $\mathcal{C}_0 = \mathcal{B}_0 = \{db\}$.

To describe small deformations of the representation~\eqref{eq:parab}, not conjugated to the initial one, we introduce a small parameter~$\delta$ and
assume that the ``deformed'' matrix for~$\eta$ is
\[
\begin{pmatrix} 1 & 1 \\ \delta & 1 \end{pmatrix}.
\]
The space $(dg)$ is thus 1-dimensional and generated by $d\delta$. A direct calculation shows that
\begin{equation}
dz_{1'} = - \left(\zeta_1^2 + \zeta_1\right)\, d\delta.
\end{equation}
So, we find
\begin{multline}
dy_{11'23} = \frac{1}{\zeta_{12}\,\zeta_{31'}} \left(\frac{dz_{11'}}{\zeta_{11'}} + \frac{dz_{23}}{\zeta_{23}} -
\frac{dz_{13}}{\zeta_{13}} - \frac{dz_{21'}}{\zeta_{21'}}\right)  \\
= - \frac{dz_{1'}}{\zeta_{31'}\,\zeta_{21'}} = \frac{\zeta_1^2 + \zeta_1}{\zeta_{31'}\,\zeta_{21'}}\, d\delta.
\end{multline}

Let us choose all the sets $\mathcal{B}_i$, $i=1, \ldots, 4$ in the same way as in the non-parabolic case. Then, an easy calculation gives
\[
\det {}_{\mathcal{B}_1} f_1 = 1,
\]
\[
\det {}_{\mathcal{B}_2} f_2 = - \frac{2}{\zeta_{12} \, \zeta_{13} \, \zeta_{23} \, \zeta_{1'2} \, \zeta_{1'3} \, \zeta_{23'}},
\]
\[
\det {}_{\mathcal{B}_3} f_3 = \zeta_{1'2} \, \zeta_{13} \, \zeta_{23},
\]
\[
\det {}_{\mathcal{B}_4} f_4 = \frac{2}{\zeta_{12} \, \zeta_{23'} \, \zeta_{1'3}}.
\]

Finally, using~\eqref{eq:tors} and~\eqref{eq:inv}, we find
\begin{equation}\label{eq:IS2S1parab}
I_\rho(S^2 \times S^1; \mathrm{PSL(2, \mathbb{C})}) = -\frac{1}{4}.
\end{equation}

\bigskip

We would like to emphasize that our answers \eqref{eq:IS2S1nonparab} and~\eqref{eq:IS2S1parab} depend on our specific choice, namely $d\lambda /
\lambda$ and~$d\delta$, of infinitesimal parameters for representation deformations.

\subsection{A relative invariant: unknots in lens spaces}
\label{subsec:relative}

In this subsection we are going to calculate our invariant for a lens space without a tubular neighborhood of unknot. Let us first briefly remind
generalities on lens spaces and their triangulations. Let $p, q$ be two coprime integers such that $0 < p < |q|$. The lens space $L(p,q)$ is defined
as the quotient manifold $S^3 / \sim$, where $\sim$ denotes the action of the cyclic group $\mathbb{Z}_p$ on $\mathbb{C}^2\supset S^3$ given by:
\[
\zeta \cdot (z_1, z_2) = (\zeta z_1, \zeta^q z_2), \quad \zeta = e^{2\pi i/p}.
\]
As a consequence the universal cover of lens spaces is the three-dimensional sphere $S^3$ and
\begin{equation}
\pi_1 \bigl( L(p,q) \bigr) \cong H_1 \bigl( L(p,q) \bigr) \cong \mathbb Z_p. \label{eq_Z_p}
\end{equation}

\begin{figure}[ht]
\begin{center}
\includegraphics[scale=0.3]{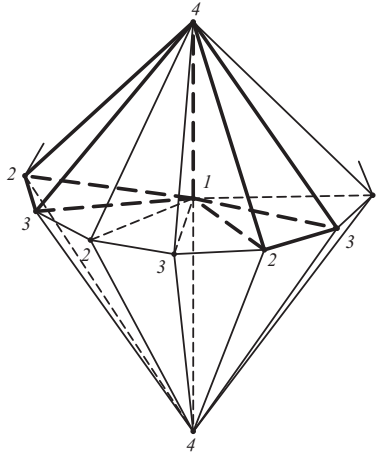}
\end{center}
\caption{A chain of two tetrahedra in a lens space} \label{fig:chain}
\end{figure}

Now we describe a triangulation of $L(p,q)$ which will be used in our calculations. Consider the bipyramid of Figure~\ref{fig:chain},
which contains $p$ vertices~$2$ and $p$ vertices~$3$. The lens space $L(p,q)$ is obtained by glueing the upper half of its surface to the lower half,
the latter having been rotated around the vertical axis through the angle $2\pi q/p$ in such way that every ``upper'' triangle $234$ is glued to some
``lower'' triangle $234$ (the vertices of the same names are identified).

A generator of the fundamental group can be represented, e.g., by some broken line $232$ (the two end points~$2$ are different) lying in
the equator of the bipyramid. We assume that a generator chosen in such way corresponds to the element $1\in \mathbb Z_p$ under the isomorphism~\eqref{eq_Z_p}.

The boldface lines (solid and dashed) in Figure~\ref{fig:chain} single out two {\em identically oriented\/} tetrahedra~$1234$ which form a chain exactly
like the one in the paper~\cite{DKM}. Going along the chain of tetrahedra in Figure~\ref{fig:chain} (e.g., along the way $212$) corresponds to the element $2\in \mathbb Z_p$ (or to $ - 2\in \mathbb Z_p$, if we go in the opposite direction). It is clear that one can also choose a pair of tetrahedra corresponding to {\em any nonzero\/} element $m$ from $H_1 \bigl( L(p,q) \bigr) \cong \mathbb Z_p$.

A knot in $L(p,q)$ determined by a tetrahedron chain of the kind of Figure~\ref{fig:chain}, i.e., going along a line like~$212$, can be called, somewhat loosely, an ``unknot'' in $L(p,q)$. It differs from any other conceivable knot, going along which gives the same element of~$H_1\bigl( L(p,q) \bigr)$, in its ``minimal knottedness'' in the following sense: the full preimage of this knot in the universal cover of space $L(p,q)$, i.e., sphere~$S^3$, being decomposed in a connected sum of simple knots, contains the {\em smallest number of summands}. Indeed, the line $212$ is equivalent, as a knot, to the segment of the straight line joining the two points~$2$; if, on the other hand, we tie a nontrivial knot on this segment, there will appear $p$ new summands in the full preimage (in the sense of connected summation) equivalent to this nontrivial knot.

\begin{figure}
\begin{center}
\includegraphics[scale=0.3]{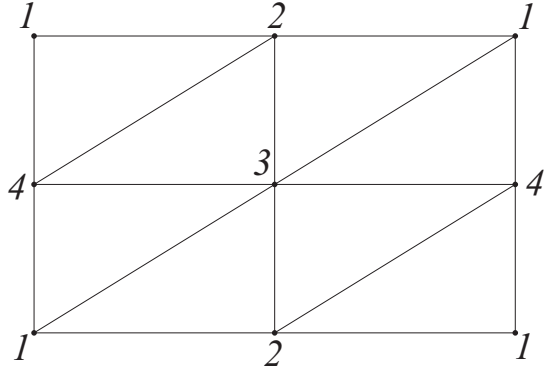}
\end{center}
\caption{Triangulation of the unknot exterior boundary} \label{8triangles}
\end{figure}

The fixed triangulation of the toric boundary of unknot exterior is presented in figure~\ref{8triangles}.
Note that the edges $14$ and~$23$ are of course ``doubled'' here, taking into account the fact that the
exterior of the unknot lies on both sides of any of them.

For the triangulation of a lens space, described above, let us consider the following algebraic complex
\begin{equation}\label{eq:complex_rel}
0 \xrightarrow{} (dy)_{\mathrm{int}} \xrightarrow{f_3} (d\varphi)_{\mathrm{int}} \oplus (d\varphi)_{\mathcal{D}} \xrightarrow{} 0.
\end{equation}
Here $(dy)_{\mathrm{int}}$ and $(d\varphi)_{\mathrm{int}}$ are subspaces of vector spaces $(dy)$ and $(d\varphi)$ corresponding to the
\emph{interior} tetrahedra and edges respectively, while $\mathcal{D}$ is a set of \textit{four} edges from the boundary and
$(d\varphi)_{\mathcal{D}}$ is a restriction of $(d\varphi)$ to this set. The number four ensures that the Euler characteristic of algebraic
complex~\eqref{eq:complex_rel} is zero. Thus, our matrix $f_3$ is square of dimension~$4p - 2$. The complex~\eqref{eq:complex_rel} is acyclic
provided $\det f_3 \neq 0$.

Note that the two distinguished tetrahedra in Figure~\ref{fig:chain} are turned into each other under a rotation through angle
$2\cdot\frac{2\pi}{p}$; similarly, an ``unknot going along the element $n\in\mathbb Z_p = H_1 \bigl(L(p,q)\bigr)$'' is determined by two tetrahedra
which differ in a rotation through angle $n\cdot\frac{2\pi}{p}$. For a different basis element in~$H_1$, this number~$n$ would change, but we are
considering the lens space~$L(p,q)$ as constructed in a fixed way from the given bipyramid in Figure~\ref{fig:chain}. We also identify $n\in \mathbb
Z_p$ with one of positive integers $1, \ldots, p-1$ (of course, $n \neq 0$). The complex~\eqref{eq:complex_rel} depends, besides the set
$\mathcal{D}$, on this number~$n$.

According to the form of~\eqref{eq:complex_rel}, the invariant comes out to be as:
\begin{equation}\label{eq:inv_rel}
I_{\mathcal{D}, n}\bigl(L(p, q); \mathrm{PSL}(2, \mathbb{C})\bigr) =  \frac{\det f_3}{\prod' \zeta^2_{ij}},
\end{equation}
where the product is taken over all edges from the triangulation not belonging to the set~$\mathcal{D}$. The value~\eqref{eq:inv_rel} remains
unchanged under any simplicial transformation of the triangulation of lens space, not involving two distinguished tetrahedra. This can be proved by
analogy with the methods of paper~\cite{DKM}.

\begin{remark}
Invariant~\eqref{eq:inv_rel} can depend \emph{a priori} on the geometry of the tetrahedron~$1234$, that is on the values $\zeta_i$ for $i = 1, 2, 3,
4$.
\end{remark}

For a given number $n$, there are in principle $\binom{12}{4} = 495$ expressions for our invariant (including zeros), depending on the set
$\mathcal{D}$. Table~\ref{table-lenses}
\begin{table}
\begin{center}
\begin{tabular}{|>{\centering}m{37.5mm}|c|c|}
\hline Set $\mathcal{D}$ & $L(7, 1)$ & $L(7, 2)$ \\\hline \vspace{1mm}
\includegraphics[scale=0.15]{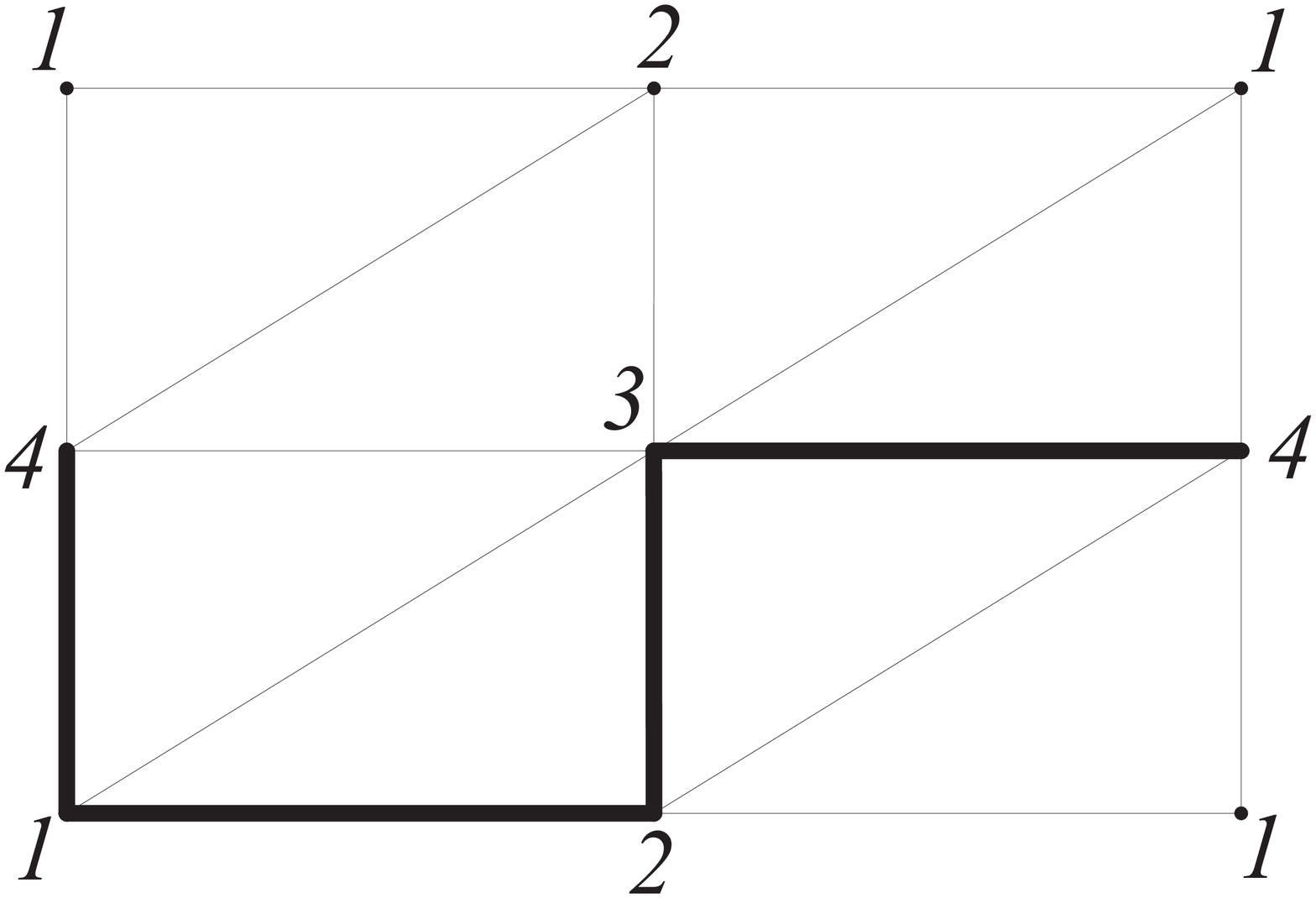} & $\displaystyle \begin{pmatrix} 6 \\ 90 \\ 300 \end{pmatrix}$ & $\displaystyle \begin{pmatrix} 48 \\ 20 \\ 6 \end{pmatrix}$ \\\hline \vspace{1mm}
\includegraphics[scale=0.15]{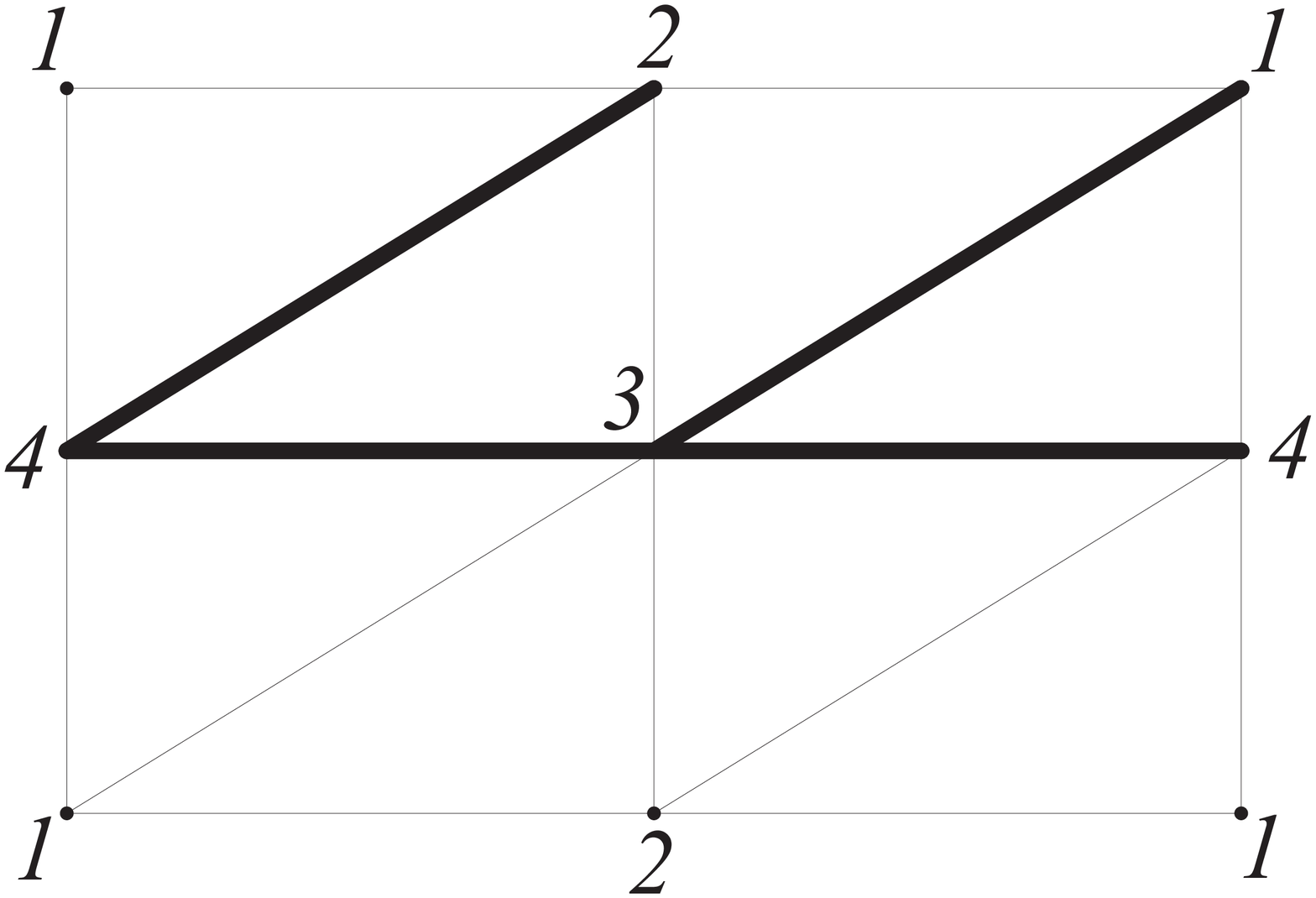} & $\displaystyle \begin{pmatrix} 49 \\ 147 \\ 245 \end{pmatrix} \frac{\zeta_{12}^2}{\zeta_{34}^2}$ & $\displaystyle \begin{pmatrix} 196 \\ 98\\ 49 \end{pmatrix} \frac{\zeta_{12}^2}{\zeta_{34}^2}$ \\\hline \vspace{1mm}
\includegraphics[scale=0.15]{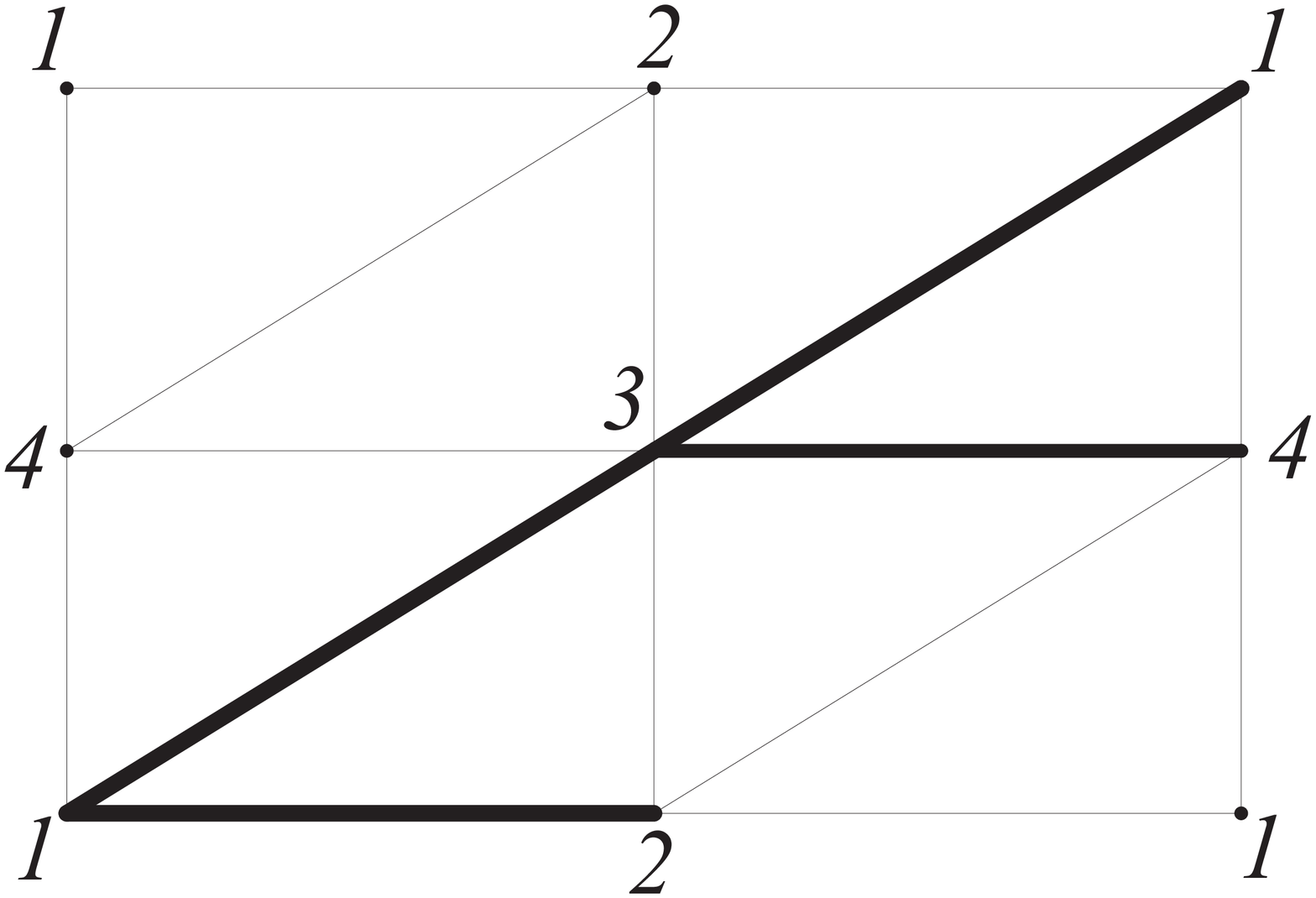} & $\displaystyle \begin{pmatrix} 294 \\ 490 \\ 588 \end{pmatrix} \frac{\zeta_{24}^2}{\zeta_{13}^2}$ & $\displaystyle \begin{pmatrix} 147 \\ 245 \\ 294 \end{pmatrix} \frac{\zeta_{24}^2}{\zeta_{13}^2}$ \\\hline \vspace{1mm}
\includegraphics[scale=0.15]{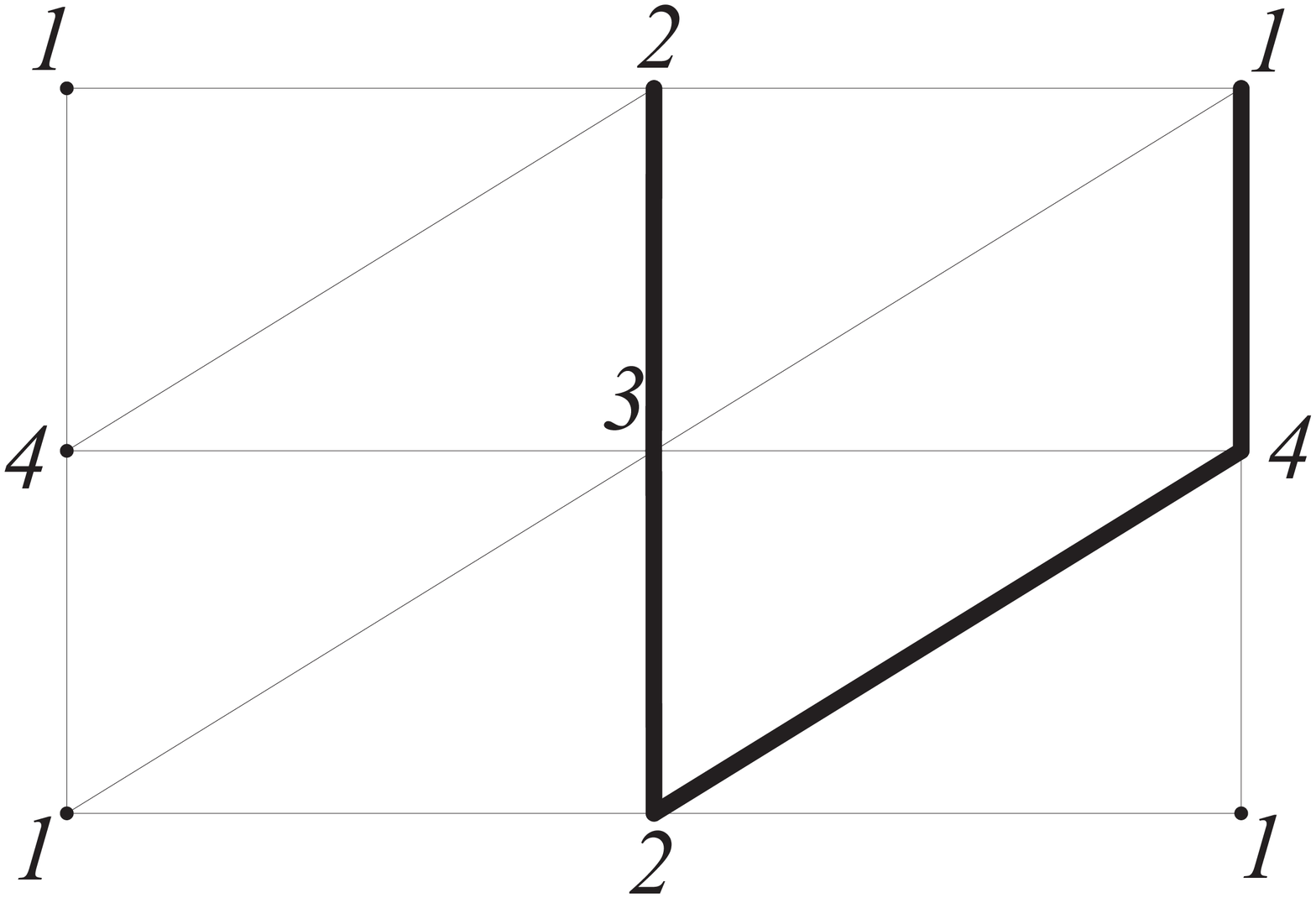} & $\displaystyle \begin{pmatrix} 1 \\ 27 \\ 125 \end{pmatrix} \frac{\zeta_{13} \zeta_{14}}{\zeta_{23} \zeta_{24}}$ & $\displaystyle \begin{pmatrix} 64 \\ 8 \\ 1 \end{pmatrix} \frac{\zeta_{13} \zeta_{14}}{\zeta_{23} \zeta_{24}}$ \\\hline \vspace{1mm}
\includegraphics[scale=0.15]{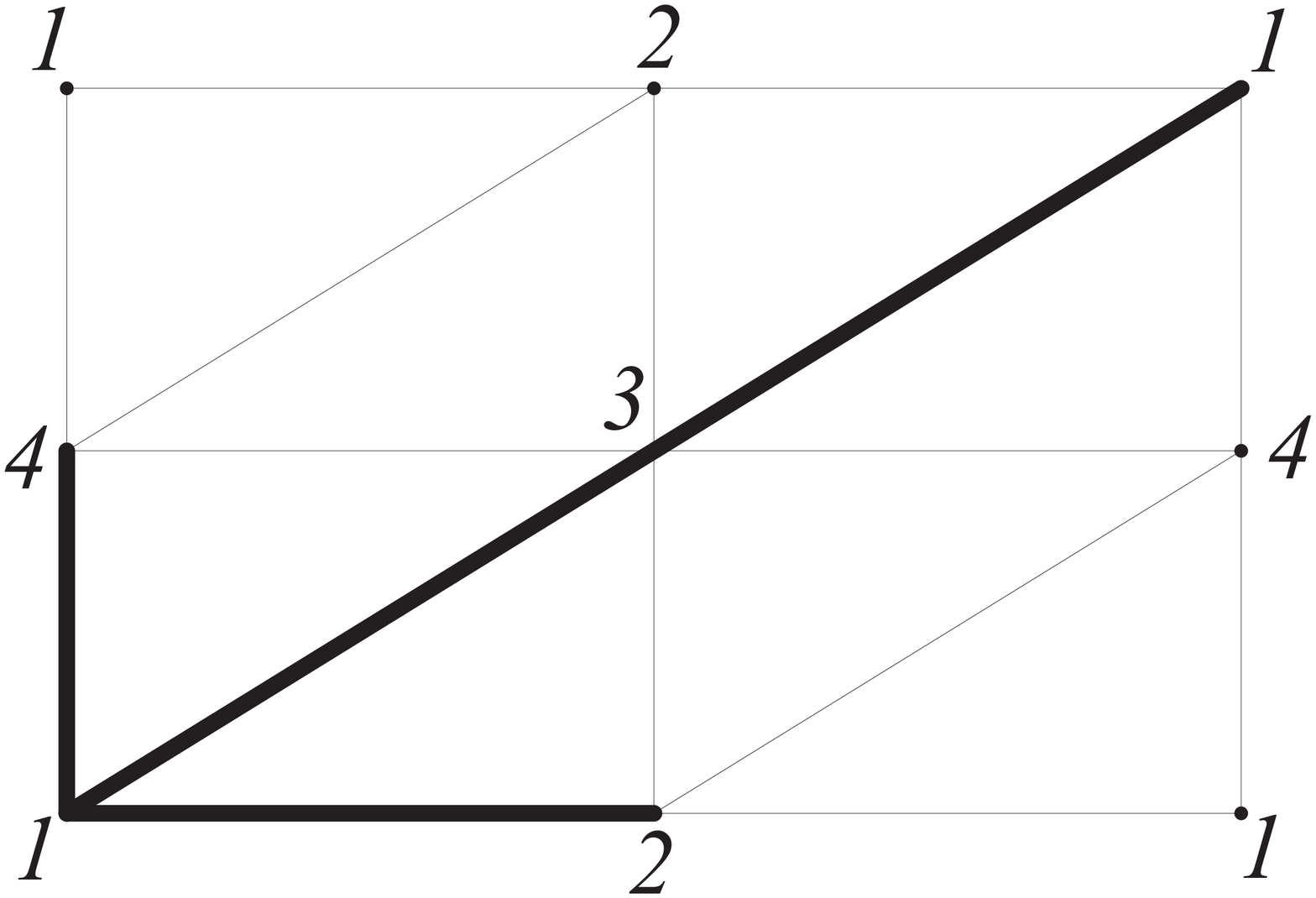} & $\displaystyle \begin{pmatrix} 42 \\ 210 \\ 420 \end{pmatrix} \frac{\zeta_{23}\zeta_{34}\zeta_{24}^2}{\zeta_{14}\zeta_{12}\zeta_{13}^2}$ & $\displaystyle \begin{pmatrix} 84 \\ 70 \\ 42 \end{pmatrix} \frac{\zeta_{23}\zeta_{34}\zeta_{24}^2}{\zeta_{14}\zeta_{12}\zeta_{13}^2}$ \\\hline
\end{tabular}
\end{center}
\caption{Some results of calculations for unknots in lens spaces} \label{table-lenses}
\end{table}
shows some results of calculation for lens spaces $L(7, 1)$ and $L(7, 2)$. The result is presented as a three-component vector $v_n = I_{\mathcal{D},
n}\bigl(L(p, q); \mathrm{PSL}(2, \mathbb{C})\bigr)$ for $n = 1, 2, 3$. The set $\mathcal{D}$ consists of boldface edges.

\section{Discussion}
\label{sec:discussion}

Here are some final remarks.
\begin{itemize}
\item Subsection~\ref{subsec:relative} shows that our construction provides a huge number of invariants, given a manifold, a knot in it, and a
triangulation of the boundary of its tubular neighborhood. Note, moreover, that a triangulation of the kind of figure~\ref{8triangles} is obviously
determined just by a \emph{framing} of the knot. Note also that, in our case, any single (component of) invariant corresponds to a subset of four
edges in figure~\ref{8triangles}, while in the ``Euclidean'' case of papers~\cite{K-tqft1,K-tqft2} it corresponded to \emph{two} subsets of equal,
but arbitrary, number of edges.
\item To produce a full-fledged topological quantum field theory like in papers~\cite{K-tqft1,K-tqft2}, we must consider how our invariants behave under a
gluing of manifolds by components of their boundaries. This part of work is left for future papers.
\item As we are using cross-ratios and ``deficit angles'' well-known in hyperbolic geometry, further research may show a deeper connection between it and
our paper.
\item The mysterious gluing of two apparently different complexes in sections \ref{sec:cr} and~\ref{right} in a single algebraic complex expects its
clarification.
\end{itemize}


\begin{thebibliography}{99}
\bibitem{Atiyah}
M. Atiyah, The geometry and physics of knots, Cambridge University Press, 1990.
\bibitem{twisted-euclidean1}
I.G. Korepanov and E.V. Martyushev, Distinguishing three-dimensional lens spaces L(7,1) and L(7,2) by means of classical pentagon equation. J.
Nonlinear Math. Phys. 9, no. 1 (2002) 86--98. arXiv:math/0210343
\bibitem{DKM}
J. Dubois, I.G. Korepanov and E.V. Martyushev, Euclidean geometric invariant of framed knots in manifolds. arXiv:math/0605164
\bibitem{K-tqft1}
I.G. Korepanov, Geometric torsions and invariants of manifolds with triangulated boundary. Accepted for publication in Theor. Math. Phys.
arXiv:0803.0123
\bibitem{K-tqft2}
I.G. Korepanov, Geometric torsions and an Atiyah-style topological field theory. Accepted for publication in Theor. Math. Phys. arXiv:0806.2514
\bibitem{SL2-1}
I.G. Korepanov and E.V. Martyushev, A classical solution of the pentagon equation related to the group SL(2). Theor. Math. Phys. 129, no. 1 (2001)
1320--1324.
\bibitem{SL2-2}
I.G. Korepanov, SL(2)-solution of the pentagon equation and invariants of three-dimensional manifolds. Theor. Math. Phys. 138, no. 1 (2004) 18--27.
arXiv:math/0304149
\bibitem{man3geo4}
I.G. Korepanov, Invariants of three-dimensional manifolds from four-dimensional Euclidean geometry. arXiv:math/0611325
\bibitem{33}
I.G. Korepanov, Euclidean 4-simplices and invariants of four-dimensional manifolds: I.~Moves $3\to 3$. Theor.Math.Phys. 131 (2002) 765--774.
arXiv:math/0211165
\bibitem{24}
I.G. Korepanov, Euclidean 4-simplices and invariants of four-dimensional manifolds: II.~An algebraic complex and moves $2\leftrightarrow 4$.
Theor.Math.Phys. 133 (2002) 1338--1347. arXiv:math/0211166
\bibitem{15}
I.G. Korepanov, Euclidean 4-simplices and invariants of four-dimensional manifolds: III.~Moves $1\leftrightarrow 5$ and related structures.
Theor.Math.Phys. 135 (2003) 601--613. arXiv:math/0211167
\bibitem{twisted-euclidean2}
E.V. Martyushev, Euclidean simplices and invariants of three-manifolds: a modification of the invariant for lens spaces. arXiv:math/0212018
\bibitem{M-thesis}
E.V. Martyushev, Geometric invariants of three-dimensional manifolds, knots and links. Ph.D. Thesis (in russian). http://www.susu.ac.ru/file/thesis.pdf
\bibitem{Pach}
U. Pachner, PL homeomorphic manifolds are equivalent by elementary shellings. Europ. J. Combinatorics 12 (1991), 129--145.
\bibitem{Tur01}
V.G. Turaev, Introduction to combinatorial torsions. Boston:~Birkh\"auser, 2000.
\end{thebibliography}
\end{document}